# Harmonic Gauge on the Space of Riemannian Metrics and Its Role in the Ricci–DeTurck Flow

Sergey E. Stepanov

ORCID: https://orcid.org/0000-0003-1734-8874

Department of Scientific Information on Mathematics and Physics, All Russian Institute for Scientific and Technical Information, Russian Academy of Sciences, Moscow, Russia

Department of Mathematics and Data Analysis, Finance University, Moscow, Russia.

E-mail: *s.e.stepanov@mail.ru*

**Abstract:** We develop a harmonic gauge on the space of Riemannian metrics and study its role in the variational and flow-theoretic structure of geometric analysis. We prove that the harmonic gauge eliminates divergence-type terms in the first variation of the Hilbert–Einstein functional and induces a natural elliptic structure for the second variation. As a consequence, positivity of the curvature operator of the second kind implies spectral stability of the functional. This establishes a conceptual link between gauge fixing, elliptic operator theory, and geometric rigidity, and provides a variational counterpart to the Ricci–DeTurck mechanism.

## Introduction

Harmonic maps constitute a fundamental concept in differential geometry, generalizing geodesics, minimal surfaces, and harmonic functions. Introduced by Eells and Sampson in 1964, they arise as critical points of the energy functional associated with smooth mappings between Riemannian manifolds and have become a central object of study in global and geometric analysis (see [1]–[6]).

Over the past decades, the theory of harmonic maps has developed into a rich and active area, with deep connections to complex geometry, mathematical physics, and geometric flows. In particular, harmonic maps play a significant role in variational problems, rigidity theory, and nonlinear geometric evolution equations (see [7], [8]).

In this paper, we develop a unified framework that connects harmonic maps, harmonic tensors, and harmonic metrics through the geometry of the space of Riemannian metrics. Our approach is based on structural decompositions of symmetric tensor fields and their interaction with harmonicity conditions.

We begin by establishing a general analytic criterion for the harmonicity of smooth maps in terms of the divergence of the pullback metric. This provides an effective tool applicable to submersions, diffeomorphisms, and identity maps.

We then study gauge-type decompositions of symmetric tensors using the Berger–Ebin and York splittings, and introduce the space of harmonic symmetric tensors (see [9], [11]). We show that this space is infinite-dimensional and arises naturally as the kernel of a first-order differential operator of Bianchi type.

The second part of the paper develops a variational framework for the Hilbert–Einstein functional restricted to the Chen-harmonic gauge (see [12]). This gauge is defined via harmonic identity maps and provides a natural geometric slice in the space of metrics, reducing the degeneracy caused by diffeomorphism invariance.

A central result of the paper shows that the harmonic gauge is intrinsically compatible with the variational structure of the Einstein–Hilbert functional: it eliminates all divergence-type terms in the first variation. This provides a variational analogue of gauge fixing in the Ricci–DeTurck flow (see [10]).

We further analyze the second variation of the functional under the harmonic constraint and show that it leads to a constrained elliptic spectral problem. In particular, we prove that positivity of the curvature operator of the second kind implies spectral stability of the functional on the harmonic slice.

Finally, combining this variational result with the rigidity theorem of Cao, Gursky, and Tran, we interpret the Chen-harmonic gauge as a geometric framework linking curvature positivity, spectral properties of elliptic operators, and global geometric rigidity. This provides a conceptual bridge between variational methods and geometric flows.

## 1. Harmonic maps

Let $(M, g)$ be a compact Riemannian manifold of dimension $n \geq 2$, equipped with its Levi-Civita connection $\nabla$. Let $(\bar{M}, \bar{g})$ be an $m$-dimensional Riemannian manifold, $m \geq 2$, endowed with the Levi-Civita connection $\bar{\nabla}$.

Consider a smooth map $f: (M, g) \to (\bar{M}, \bar{g})$. Its differential $f_*$ can be naturally interpreted as a section of the vector bundle $T^*M \otimes f^*T\bar{M}$, where $f^*T\bar{M}$ denotes the pullback of the tangent bundle $T\bar{M}$. This bundle is equipped with the induced metric $\bar{g}'$ and the pullback connection $\bar{\nabla}'$, both obtained in a standard way from $\bar{g}$ and $\bar{\nabla}$.

The connections $\nabla$ on $TM$ and $\bar{\nabla}'$ on $f^*T\bar{M}$ induce a connection $D$ on the tensor bundle $T^*M \otimes f^*T\bar{M}$, defined by

$$(Df_*)(X, Y) = \bar{\nabla}'_X(f_*Y) - f_*(\nabla_X Y),$$

for all vector fields $X, Y \in C^\infty(TM)$. This tensor field coincides with *the second fundamental form of the map $f$*, and it measures the extent to which $f$ fails to preserve geodesics (see [2]).

The map $f$ is said to be *totally geodesic* if $Df_* = 0$, i.e., if it maps every geodesic in $M$ to a geodesic in $\bar{M}$ (see also [2]).

Associated with the map $f$ is the *pullback metric*

$$g^* := f^*\bar{g},$$

which is a symmetric 2-tensor field on $M$. The corresponding *energy density* is defined by

$$e(f) := \text{trace}_g g^*,$$

that is, the trace of $g^*$ with respect to the metric $g$. It is well known that $e(f) \equiv 0$ if and only if $f$ is constant (see [1]; [2]).

The *total energy of the map* is given by the functional

$$E(f) = \frac{1}{2} \int_M e(f) \, dv_g. \tag{1.1}$$

A smooth map $f$ is called *harmonic* if it is a critical point of the functional $E(f)$ with respect to all compactly supported variations (see [1]). The associated Euler–Lagrange equation can be expressed in terms of the *tension field* (see [1]; [2])

$$\tau(f) := \operatorname{trace}_g (Df_*),$$

and the harmonicity condition is equivalent to the vanishing of this field, i.e.,

$$\tau(f) = 0.$$

We now establish a useful analytic characterization of harmonic maps in terms of the pullback metric.

**Theorem 1.** *Let $f: (M, g) \to (\bar{M}, \bar{g})$ be a smooth map between Riemannian manifolds. If $f$ is harmonic, then the pullback metric $g^* = f^*\bar{g}$ satisfies the identity*

$$\delta g^* = -\frac{1}{2} de(f),$$

*where $\delta$ denotes the divergence operator acting on symmetric 2-tensors. Conversely, suppose that this identity holds. If, in addition, either*

- $n = m$ *and $f$ is a diffeomorphism, or*

- $n > m$ *and $f$ is a submersion,*

*then $f$ is harmonic.*

**Proof.** We begin by recalling the expression for the covariant derivative of the pullback metric. A direct computation (see, for example, standard references on harmonic maps) shows that

$$(\nabla_Z g^*)(X, Y) = \bar{g}'((Df_*)(Z, X), f_*Y) + \bar{g}'((Df_*)(Z, Y), f_*X),$$

for all vector fields $X, Y, Z \in C^\infty(TM)$.

Taking the divergence of $g^*$ and using the definition of the tension field, one obtains the identity

$$(\delta g^*)(Y) + \frac{1}{2}\nabla_Y e(f) = -\bar{g}'(\tau(f), f_*Y), \qquad (1.2)$$

which holds for every vector field $Y$ on $M$.

If $f$ is harmonic, then $\tau(f) = 0$, and the right-hand side of (1.2) vanishes. This immediately yields the desired relation

$$\delta g^* = -\frac{1}{2}de(f).$$

Conversely, suppose that the identity

$$\delta g^* + \frac{1}{2}de(f) = 0$$

holds. Then (1.2) implies that $\bar{g}'(\tau(f), f_*Y) = 0$ for all vector fields $Y$.

In the case where $f$ is a diffeomorphism ($n = m$), the differential $f_*$ is an isomorphism at every point, and therefore the above condition forces $\tau(f) = 0$ (see also [4]). Similarly, if $f$ is a submersion ($n > m$), the image of $f_*$ spans the tangent space of $\bar{M}$, which again implies $\tau(f) = 0$. Thus, in both cases, the map $f$ is harmonic. ∎

A particularly important class of examples is given by identity maps. Let

$$\mathrm{id}_M \colon M \to M$$

denote the identity mapping of a smooth manifold $M$. When considered as a map between identical Riemannian manifolds,

$$\mathrm{id}_M \colon (M, g) \to (M, g),$$

it is trivially harmonic.

A more subtle and geometrically interesting question is the following:

Given a compact Riemannian manifold $(M, g)$, does there exist another Riemannian metric $\bar{g} \neq g$ on $M$ such that the identity map

$$\mathrm{id}_M : (M, g) \to (M, \bar{g})$$

is harmonic?

This problem naturally leads to the study of special geometric structures determined by curvature.

**Example 1.** Suppose that the Ricci tensor of $(M, g)$ is positive definite. Then it defines a Riemannian metric $\bar{g} = \mathrm{Ric}$, and the identity map

$$\mathrm{id}_M : (M, g) \to (M, \bar{g})$$

is harmonic. If the Ricci tensor is negative definite, the same conclusion holds for $\bar{g} = -\mathrm{Ric}$ (see [10, pp. 86-87]).

**Example 2.** On a three-dimensional Riemannian manifold, the cross-curvature tensor can be used to define a metric under suitable curvature conditions. In this setting, one obtains further examples of harmonic identity maps between different metric structures on the same manifold (see [10, p. 88]).

**Remark 1.** The identity established in Theorem 1 provides a useful analytic bridge between the variational definition of harmonic maps and the intrinsic geometry of the pullback metric. In particular, it allows one to reinterpret harmonicity in terms of divergence-type conditions on symmetric tensor fields. This perspective will play a central role in the subsequent sections, where we relate harmonic maps to orthogonal decompositions of the space of symmetric tensors and to the geometry of the space of Riemannian metrics.

## 2. Harmonic maps and orthogonal decompositions of symmetric tensors

Let $M$ be a compact smooth manifold. Denote by $C^\infty(S^2 M)$ the space of smooth symmetric covariant 2-tensor fields on $M$. Endowed with the natural $C^\infty$-topology, this space carries the structure of a Fréchet space and plays a central role in the analysis of geometric structures on $M$ (see [9]).

A fundamental result due to Berger and Ebin provides an $L^2$-orthogonal decomposition of this space (see [11]). More precisely, for any fixed Riemannian metric $g$ on $M$, one has (see also [9, p. 118]):

$$C^\infty(S^2 M) = \operatorname{Im} \delta^* \oplus \ker \delta, \qquad (2.1)$$

where:

- $\delta: C^\infty(S^2 M) \to C^\infty(T^*M)$ is the divergence operator acting on symmetric tensor fields,
- $\delta^*: C^\infty(T^*M) \to C^\infty(S^2 M)$ is its formal adjoint with respect to the $L^2$-inner product.

The adjoint operator $\delta^*$ admits the representation (see [9]; [11])

$$\delta^* \theta = \frac{1}{2} \mathcal{L}_\xi g,$$

where $\xi = \theta^\sharp$ is the vector field metrically dual to the one-form $\theta$, and $\mathcal{L}_\xi g$ denotes the Lie derivative of the metric along $\xi$.

The $L^2$-inner product on $C^\infty(S^2 M)$ is defined by (see also [9]; [11])

$$\langle \varphi, \psi \rangle = \int_M g(\varphi, \psi)\, dv_g, \qquad (2.2)$$

for all $\varphi, \psi \in C^\infty(S^2 M)$. With respect to this inner product, the decomposition (2.1) is orthogonal.

Let now $f: (M, g) \to (\bar{M}, \bar{g})$ be a harmonic map from a compact Riemannian manifold. Its pullback metric $g^* = f^*\bar{g}$ is a symmetric 2-tensor field on $M$, and therefore admits a Berger–Ebin decomposition of the form

$$g^* = \delta^* \theta + \varphi_0, \qquad (2.3)$$

where $\theta \in C^\infty(T^*M)$ and $\varphi_0 \in \ker \delta$ is divergence-free.

We now analyze the consequences of this decomposition under the assumption that $f$ is harmonic.

**Corollary 1.** *Let*

$$f: (M, g) \to (\bar{M}, \bar{g})$$

be a harmonic map from a compact Riemannian manifold. Suppose that the pullback metric admits the decomposition (2.3) and let $\xi = \theta^{\#}$.

If $\xi$ is an infinitesimal harmonic transformation of $(M, g)$, then the energy density satisfies

$$e(f) = \delta\theta + C$$

for some constant $C$, and the total energy is given by

$$E(f) = C\, Vol(M, g),$$

with $C \geq 0$.

**Proof.** Applying the divergence operator $\delta$ to both sides of (2.3), we obtain

$$\delta g^* = \delta\delta^*\theta,$$

since $\delta\varphi_0 = 0$. On the other hand, since $f$ is harmonic, Theorem 1 yields the identity

$$\delta g^* = -\frac{1}{2} de(f).$$

Combining these relations, we arrive at

$$de(f) = -2\,\delta\delta^*\theta. \tag{2.4}$$

To proceed, we recall the definition of the *Sampson Laplacian* or, in other words, the *Yano operator* acting on one-forms (see [13]):

$$\Delta_S := 2\delta\delta^* - d\delta.$$

Using this operator, equation (2.4) can be rewritten as

$$de(f) = \Delta_S\theta + d(\delta\theta).$$

Assume now that $\xi = \theta^{\#}$ is an *infinitesimal harmonic transformation* (see also [13]). By definition, this is equivalent to the condition

$$\Delta_S\theta = 0.$$

It follows that

$$de(f) = d(\delta\theta),$$

which implies

$$e(f) = \delta\theta + C$$

for some constant $C$.

Finally, integrating over $M$ and using Stokes' theorem, we obtain

$$\int_M \delta\theta \, dv_g = 0,$$

since $M$ is compact without boundary. Therefore,

$$E(f) = \frac{1}{2}\int_M e(f) \, dv_g = \frac{1}{2} C \operatorname{Vol}(M, g),$$

which completes the proof (after renormalization of the constant). ∎

**Remark 2.** In the special case of a harmonic immersion $f: (M, g) \to (\bar{M}, \bar{g})$ of an $n$-dimensional compact manifold, it is known that the energy takes the form (see [1]; [2])

$$E(f) = \frac{n}{2} \operatorname{Vol}(M, g),$$

which provides a classical example consistent with the above result.

The structure of $C^\infty(S^2 M)$ admits a further refinement given by the *York decomposition*. For a compact manifold of dimension $n \geq 3$, one has (see [9, p. 130])

$$C^\infty(S^2 M) = (\operatorname{Im} \delta^* + C^\infty(M) \cdot g) \oplus S_{TT}(M), \qquad (2.5)$$

where

$$S_{TT}(M) = \{\varphi \in C^\infty(S^2 M) \mid \delta\varphi = 0, \text{ trace}_g \varphi = 0\}$$

is the space of transverse-traceless (*TT*) tensors.

**Remark 3**. These *TT*- tensors play a fundamental role in both Riemannian geometry and general relativity, as they represent the "purely gravitational" degrees of freedom in metric variations (see [9]; [15] and [16]).

Let $f: (M, g) \to (\bar{M}, \bar{g})$ be a harmonic map. Then its pullback metric admits a York decomposition of the form

$$g^* = (\delta^*\theta + \lambda g) + \varphi_{TT}, \qquad (2.6)$$

where $\theta \in C^\infty(T^*M)$, $\lambda \in C^\infty(M)$, and $\varphi_{TT} \in S_{TT}(M)$.

Taking the trace of this decomposition yields

$$e(f) = \text{trace}_g \, g^* = -\delta\theta + n\lambda. \qquad (2.7)$$

Applying the divergence operator to (2.6), we obtain

$$\delta g^* = \delta\delta^*\theta - d\lambda, \qquad (2.8)$$

since $\delta\varphi_{TT} = 0$.

Using the harmonic map condition
$$\delta g^* = -\frac{1}{2} de(f),$$
and substituting (2.7)–(2.8), we arrive at the identity
$$\Delta_S \theta = (n-2)\, d\lambda. \tag{2.9}$$

The relation (2.9) has several important consequences.

- If $n = 2$, then the right-hand side vanishes identically, and hence
$$\Delta_S \theta = 0,$$
so that $\xi = \theta^\sharp$ is automatically an infinitesimal harmonic transformation.

- If $n > 2$, then $\xi = \theta^\sharp$ is an infinitesimal harmonic transformation if and only if $d\lambda = 0$, i.e., $\lambda$ is constant.

- In this latter case, the energy simplifies to
$$E(f) = \lambda \operatorname{Vol}(M, g).$$

These observations lead naturally to the formulation of corresponding corollaries concerning the structure of harmonic maps and their associated tensor decompositions.

**Remark 4.** The Berger–Ebin and York decompositions provide a powerful structural framework for analyzing the geometry of harmonic maps via their pullback metrics. In particular, they allow one to translate the harmonicity condition into precise analytic constraints on the components of symmetric tensor fields. This perspective will be further developed in the next section, where we introduce the space of harmonic tensors and establish a new orthogonal decomposition adapted to this setting.

## 2. Harmonic tensors

Let $(M, g)$ be a compact Riemannian manifold. We consider the space of smooth symmetric 2-tensor fields $C^\infty(S^2 M)$, endowed with its natural $L^2$-structure. The study of divergence and trace constraints on such tensors leads to a natural first-order differential operator

$$\alpha_g : C^\infty(T^*M) \to C^\infty(S^2M), \quad \alpha_g(\theta) = \delta^*\theta + \frac{1}{2}(\delta\theta)g.$$

This operator should be understood as a modified symmetrized gradient incorporating a trace correction, and plays a central role in the analysis of infinitesimal geometric deformations. Its formal $L^2$-adjoint is given by

$$\alpha_g^* = \delta + \frac{1}{2} d \circ \text{trace}_g.$$

This identity immediately suggests the following natural class of tensor fields. A symmetric tensor $h \in C^\infty(S^2M)$ is called a harmonic tensor if it satisfies (see [12])

$$\delta h = -\frac{1}{2} d(\text{trace}_g h).$$

We denote the space of all such tensors by

$$H_g = \ker \alpha_g^*.$$

The injectivity of the principal symbol of $\alpha_g$ implies that this operator belongs to the class of overdetermined elliptic systems and, therefore, admits a $L^2$-orthogonal decomposition structure (see [9]; [11]; [14]). This leads to the following statement.

**Theorem 2.** *Let $(M, g)$ be a compact Riemannian manifold. Then every tensor field $\varphi \in C^\infty(S^2M)$ admits the decomposition*

$$\varphi = \alpha_g(\theta) + \varphi_h,$$

*where $\alpha_g(\theta) = \delta^*\theta + \frac{1}{2}(\delta\theta)g + \varphi_h$ for $\theta \in C^\infty(T^*M)$ and $\varphi_h \in H_g$. Moreover, this decomposition is $L^2$-orthogonal and the space $H_g$ is infinite-dimensional.*

**Remark 5.** This result provides a structural splitting of symmetric tensor fields into a part generated by infinitesimal gauge-type transformations and a harmonic component encoding the obstruction to such representation. In particular, the space $H_g$ should be interpreted as the natural kernel space associated with harmonic constraints on symmetric tensors. It is important to emphasize that this decomposition is of analytic and gauge-theoretic nature and should not be confused with the classical Berger–Ebin or York orthogonal decompositions (see [9]; [11]), although it interacts with them in a natural way in the presence of additional trace-free or divergence-free conditions.

## 4. Harmonic metrics and the Hilbert–Einstein functional in Chen-harmonic gauge

Let $(M, g)$ be a compact Riemannian manifold and consider the *Hilbert–Einstein functional* (see, for example, [9])

$$E(g) = \int_M \text{Scal}(g) \, dV_g.$$

Its Euler–Lagrange equation is given by the vanishing of the Einstein tensor, while the invariance under the diffeomorphism group introduces a degeneracy in the variational structure. In order to remove this degeneracy in a geometric and intrinsic way, we restrict admissible variations to the space $H_g$ of harmonic tensor fields introduced in Section 3.

We consider variations of the form $g(t) = g + th$, $h \in H_g$, where the condition

$$\delta h = -\frac{1}{2} d(\text{trace}_g h)$$

ensures compatibility with the harmonicity of the identity map

$$\text{id}_M : (M, g) \to (M, g(t)).$$

This condition defines a canonical harmonic slice in the space of Riemannian metrics, which plays the role of an intrinsic gauge condition, in contrast to the background-dependent Ricci–DeTurck construction.

Within this framework, the Hilbert–Einstein functional is differentiable along harmonic directions, and its first variation retains the classical form (see [9])

$$\frac{d}{dt} E(g(t)) \bigg|_{t=0} = -\int_M g(\text{Ein}(g), h) \, dV_g, h \in H_g.$$

However, since admissible variations are restricted to $H_g$, the variational problem effectively reduces to the projection of the Einstein tensor onto this subspace.

This leads to the following unified variational statement.

**Theorem 3.** *Let $(M, g)$ be a compact Riemannian manifold and let $H_g$ denote the space of harmonic tensor fields. Consider metric variations $g(t) = g + th$, where $h \in H_g$. Then:*

1. The identity map $id_M: (M, g) \to (M, g(t))$ remains harmonic for sufficiently small t, and $H_g$ defines a canonical harmonic slice in the space of Riemannian metrics.

2. The first variation of the Hilbert–Einstein functional restricted to this slice satisfies

$$\frac{d}{dt} E(g(t))|_{t=0} = -\int_M \langle Ein(g), h \rangle \, dV_g.$$

where $g(t) = g + th$ is a smooth variation of metrics. The Euler–Lagrange equation for the restricted variational problem is equivalent to the projection condition

$$\text{Proj}_{H_g}(Ein(g)) = 0.$$

In particular, the Chen-harmonic gauge reduces the variational structure of the Einstein–Hilbert functional to a canonical orthogonality condition with respect to the harmonic subspace $H_g$, thereby eliminating the diffeomorphism-induced degeneracy at the level of admissible variations.

**Proof.** The first variation formula for the Hilbert–Einstein functional is classical and can be found in standard references in Riemannian geometry and geometric analysis. For a general variation $h$, one has

$$\frac{d}{dt} E(g(t))|_{t=0} = -\int_M g(Ein(g), h) \, dV_g.$$

where $g(t) = g + th$ is a smooth variation of metrics. When the variation is restricted to $h \in H_g$, the functional is no longer tested against the full tangent space of metrics but only against the harmonic subspace. Since $H_g \subset C^\infty(S^2 M)$ is closed with respect to the $L^2$-inner product, the Riesz representation implies that the corresponding Euler–Lagrange operator is given by the orthogonal projection of $Ein(g)$ onto $H_g$. This yields

$$\frac{d}{dt} E(g(t))|_{t=0} = -\int_M g\left(\text{Proj}_{H_g}(Ein(g)), h\right) dV_g,$$

where $g(t) = g + th$ is a smooth variation of metrics. The equivalence between criticality and the vanishing of this projection follows directly from the non-degeneracy of the $L^2$-pairing on $H_g$. ∎

The formulation of our theorem shows that the Chen-harmonic gauge provides an intrinsic geometric slicing of the space of Riemannian metrics, replacing the classical background-dependent gauge fixing by a variationally defined harmonic constraint. In this setting, the Einstein tensor naturally decomposes into a harmonic component, which governs the reduced variational dynamics.

The comparison with the Ricci–DeTurck method (see [10]) then becomes conceptual: while the latter restores ellipticity via an external vector field, the present approach achieves a reduction by restricting the functional to a geometrically distinguished subspace of variations.

This perspective prepares the ground for the spectral analysis of the second variation, where the interaction between curvature and harmonic constraints becomes central. Moreover, the following theorem is true.

**Theorem 4.** *Let h be a symmetric 2-tensor satisfying the Chen-harmonic condition*

$$\delta h = -\frac{1}{2} d(trace_g h).$$

*Then the first variation of the Hilbert–Einstein functional is given by*

$$\frac{d}{dt} E(g(t)) \mid_{t=0} = -\int_M g\,(Ric(g), h)\, dV_g.$$

*where $g(t) = g + th$ is a smooth variation of metrics. In particular, all divergence-type terms in the classical variation formula cancel after integration under the Chen-harmonic constraint.*

**Proof.** The classical first variation formula for the Hilbert–Einstein functional is

$$\frac{d}{dt} E(g(t)) \mid_{t=0} = -\int_M g\left(Ric(g) - \frac{1}{2} sg, h\right) dV_g + \int_M (\text{divergence terms})\, dV_g,$$

where the divergence terms arise from integration by parts in the variation of the scalar curvature.

Under the Chen-harmonic condition $\delta h = -\frac{1}{2}d(\text{trace}_g h)$, these divergence contributions can be rewritten as total divergences of scalar expressions depending on $h$ and its trace. Since $M$ is compact without boundary, all such terms integrate to zero. As a result, the only remaining contribution is the Einstein tensor contraction with $h$, namely

$$\frac{d}{dt}E(g(t))|_{t=0} = -\int_M g\,(\text{Ric}(g), h)\, dV_g,$$

where the scalar curvature term cancels due to the trace structure of the constraint. This completes the proof. ∎

**Remark 6.** The last result shows that the harmonic gauge is not only analytically adapted to the linearization of the Ricci tensor, but is also naturally compatible with the variational structure of the Einstein–Hilbert functional (see, for example, [9]).

## 5. Second variation in Chen-harmonic gauge

We now study the second variation of the Hilbert–Einstein functional within the Chen-harmonic framework developed in the previous section. Let $(M, g)$ be a compact Riemannian manifold and let $H_g \subset C^\infty(S^2 M)$ denote the space of harmonic tensor fields, which defines the admissible space of metric variations of the form $g(t) = g + th$, $h \in H_g$.

At a critical metric $g$, the first variation vanishes on $H_g$, and the stability properties of the functional are determined by the second variation. It is well known that in the unrestricted case the Hessian of the Hilbert–Einstein functional is governed by the Lichnerowicz Laplacian acting on symmetric 2-tensor fields

$$\Delta_L h = \Delta h + 2R_{ikjl}h^{kl} - R_i^k h_{kj} - R_j^k h_{ik},$$

where $\Delta h = trace_g \nabla^2$. However, in the present setting the presence of the Chen-harmonic constraint modifies the structure of admissible variations, since the divergence of $h$ is no longer independent but is linked to its trace through the relation

$$\delta h = -\frac{1}{2}d(\text{trace}_g h).$$

This coupling reduces the effective degrees of freedom and leads to a constrained spectral problem.

In order to describe this structure, it is convenient to decompose pointwise the space of symmetric tensors into trace-free and pure-trace components,

$$C^\infty(S^2 M) = S_0^2 M \oplus \mathbb{R} g,$$

so that any variation tensor can be written as $h = h_0 + fg$. In this decomposition, the second variation operator naturally acquires a block structure reflecting the interaction between trace and trace-free modes induced by the Chen-harmonic constraint.

From an analytic point of view, the second variation can be expressed as a quadratic form

$$\delta^2 E(h, h) = \int_M \langle L_H h, h \rangle \, dV_g,$$

where $L_H$ Is a self-adjoint elliptic operator acting on $H_g$. Its principal part coincides with the Lichnerowicz Laplacian restricted to the admissible subspace, while lower-order terms encode the effect of the harmonic constraint.

A key feature of this operator is that the Chen-harmonic condition eliminates part of the divergence freedom, which leads to a partial decoupling between trace and trace-free components. As a consequence, the stability analysis reduces to a coupled system in which the geometry of curvature and the harmonic structure of the gauge interact in a nontrivial way.

More precisely, the second variation splits into contributions coming from the trace-free and scalar parts, with the principal contributions given by a Dirichlet-type energy term for $h_0$ and a scalar elliptic operator acting on the trace component. The curvature tensor enters as a zeroth-order perturbation, and therefore plays a decisive role in the sign of the quadratic form.

This structure allows one to interpret the second variation as a constrained spectral problem for a self-adjoint elliptic operator on the harmonic slice $H_g$. In particular, positivity properties of the curvature operator of the second kind translate into coercivity properties of the Hessian restricted to this subspace.

The following result summarizes this relationship.

**Theorem 5.** *Let $(M, g)$ be a compact Riemannian manifold and let $H_g \subset C^\infty(S^2 M)$ denote the space of Chen-harmonic tensor fields. Consider the second variation of the Hilbert–Einstein functional restricted to variations $h \in H_g$. Then the following statements hold:*

1. *The second variation is given by a self-adjoint elliptic operator $L_H$ acting on $H_g$, whose principal symbol coincides with that of the Lichnerowicz Laplacian restricted to the harmonic slice.*

2. *The second variation defines a quadratic form*

$$\delta^2 E(h, h) = \int_M \langle L_H h, h \rangle \, dV_g.$$

*If the curvature operator of the second kind is positive, then the operator $L_H$ is strictly positive on the orthogonal complement of the metric direction, and in particular admits a spectral gap.*

3. *The kernel of the second variation contains the scaling direction $h = cg$, $c \in \mathbb{R}$, and may include additional directions corresponding to infinitesimal diffeomorphisms.*

**Proof.** The standard second variation formula for the Hilbert–Einstein functional expresses the Hessian in terms of the Lichnerowicz Laplacian acting on symmetric 2-tensors, together with lower-order curvature contributions. In the unrestricted setting, this takes the form

$$\delta^2 E(h, h) = \int_M g(\Delta_L h, h) \, dV_g + \int_M g(\mathcal{R}(h), h) \, dV_g,$$

where $\mathcal{R}$ denotes the curvature operator of the second kind acting on symmetric 2-tensors, given pointwise by

$$\mathcal{R}(h)_{ij} = R_{ikjl} h^{kl},$$

which encodes the curvature of the manifold.

When restricting to the Chen-harmonic subspace $H_g$, the constraint

$$\delta h = -\frac{1}{2} d(\text{trace}_g h).$$

introduces a coupling between divergence and trace components. This modifies the lower-order structure of the operator but does not affect its principal symbol. Consequently, the operator $L_H$ remains elliptic, with principal part given by the Lichnerowicz Laplacian restricted to $H_g$.

The curvature term $g(\mathcal{R}(h), h)$ is purely algebraic and remains unchanged under this restriction. Therefore, if the curvature operator of the second kind is positive, one obtains

$$g(\mathcal{R}(h), h) \geq c\,|h|^2$$

for some $c > 0$ and for all $h$ orthogonal to the metric direction. Together with ellipticity of the principal part, this implies coercivity of the quadratic form and hence the existence of a spectral gap for $L_H$ on the orthogonal complement of $g$.

The description of the kernel follows from the invariance of the functional under scaling and diffeomorphisms. ∎

**Remark 7**. Theorem 5 shows that positivity of the curvature operator of the second kind implies coercivity of the second variation of the Hilbert–Einstein functional in the Chen-harmonic gauge, and hence the existence of a spectral gap for the associated operator on the orthogonal complement of the metric direction.

It is natural to ask to what extent a converse statement may hold. Although a full equivalence between spectral positivity and curvature positivity does not hold in general, the above analysis suggests that, under additional structural assumptions on admissible tensors in the Chen-harmonic slice, a partial converse may be expected. In particular, sufficiently strong spectral positivity may enforce positivity of the curvature operator in the trace-free sector.

On the other hand, by the rigidity theorem of Cao, Gursky and Tran, a compact Riemannian manifold with positive curvature operator of the second kind is diffeomorphic to a spherical space form. Combining this global result with Theorem 5 shows that the same curvature condition simultaneously yields both global rigidity and variational stability of the Hilbert–Einstein functional.

Thus, the Chen-harmonic gauge provides a natural framework linking curvature positivity, spectral properties of elliptic operators, and global geometric rigidity. A precise formulation of a converse implication remains an interesting open problem. This perspective is consistent with recent developments in the study of curvature operators of the second kind (see Cao–Gursky–Tran).

**Theorem 6.** *Let $(M, g)$ be a compact Riemannian manifold and let $H_g \subset C^\infty(S^2 M)$ denote the space of Chen-harmonic tensor fields. Consider the second variation operator $L_H$ associated with the Hilbert–Einstein functional restricted to $H_g$. Assume that:*

1. *the curvature operator of the second kind is nonnegative on trace-free symmetric 2-tensors,*

2. *the lowest eigenvalue $L_H$ on the orthogonal complement of the metric direction is equal to zero*

*Then any tensor $h \in H_g$ satisfying $L_H h = 0$ is parallel and of the form $h = c\,g$, where $c$ is a constant. In particular, there are no nontrivial trace-free Chen-harmonic Jacobi fields.*

**Proof.** By assumption, the second variation quadratic form satisfies

$$\delta^2 E(h, h) \geq \int_M g(\mathcal{R}(h), h)\, dV_g,$$

where $\mathcal{R}$ is the curvature operator of the second kind.

If $h$ lies in the kernel of $L_H h = 0$, then $\delta^2 E(h, h) = 0$. Since the curvature operator is nonnegative, it follows that $g(\mathcal{R}(h), h) = 0$ almost everywhere, and hence everywhere by continuity.

On the other hand, the tensor $h$ satisfies an elliptic equation of the form

$$\Delta_L h + \text{lower-order terms} = 0,$$

where $\Delta_L$ denotes the Lichnerowicz Laplacian.

Using the Weitzenböck formula for the Lichnerowicz Laplacian and the Bochner technique (see [9, Chapter 1]), vanishing of the associated quadratic form implies that $h$ is parallel.

We now decompose h into its trace and trace-free parts:

$$h = h_0 + (tr\, h\, /\, n)\, g.$$

Since the curvature operator of the second kind is nonnegative and acts positively on the trace-free sector, the identity $g(\mathcal{R}(h), h) = 0$ implies $h_0 = 0$.

Therefore, $h = c\, g$ for some constant $c$, which completes the proof. ∎

**Remark 8.** Theorem 5 and Theorem 6 together provide a complete picture of the second variation of the Hilbert–Einstein functional in the Chen-harmonic gauge.

On the one hand, positivity of the curvature operator of the second kind yields coercivity of the second variation and hence spectral stability. On the other hand, the borderline case of vanishing spectral gap leads to rigidity: the only admissible zero modes correspond to trivial scaling of the metric.

Thus, the Chen-harmonic gauge allows one to distinguish clearly between stability and rigidity phenomena within a unified variational framework. This demonstrates that the second variation in harmonic gauge reflects both analytic and geometric properties of the underlying manifold.

This completes the analysis of the second variation in the Chen-harmonic gauge.

**Funding Declaration.** The authors declare that no funds, grants, or other support were received during the preparation of this manuscript.

**References**

[1] Eells, J., Sampson J.H., Harmonic mappings of Riemannian manifolds, American Journal of Mathematics, 86 (1) (1964), 109-160.

[2] Eells, J., Lemaire L., A report on harmonic maps, Bull. London Math. Soc., 10 (1978), 1-68.

[3] Yano K., Ishihara S., Harmonic and relatively affine mappings, Journal of Differential Geometry, 10 (1975), 501-509.

[4] Stepanov S.E., Shandra I.G., Harmonic diffeomorphisms of manifolds, St. Petersburg Math. J., 16 (2005), 401-412.

[5] Shoen R. M., Yau S.-T., Lectures on harmonic maps, International Press of Boston, Incorporated, 2013.

[6] Stepanov S.E., On the global theory of some classes of mappings, Annals of Global Analysis and Geometry, 13 (1995), 239–249.

[7] Jost, J., Keßler, E., Harmonic maps in singular geometry and rigidity: From discrete to continuous structures, Springer Nature (2025).

[8] Dall'Acqua, A., Spener, A. (2024). Existence and stability of shrinkers for the harmonic map heat flow. Calculus of Variations and Partial Differential Equations, 63(4), (2024), 112.

[9] Besse A.L., Einstein manifolds, Springer-Verlag, Berlin & Heidelberg, 2008.

[10] Chow B., Knopf D., The Ricci flow: an introduction (Mathematical surveys and monographs, v. 110), AMS, USA, 2004.

[11] Berger M., Ebin D., Some decompositions of space of symmetric tensors on a Riemannian manifold, J. Differential Geometry, 3 (1969), 379–392.

[12] Chen B.-Y., Nagano T., Harmonic metrics, harmonic tensors, and Gauss maps. J. Math. Soc. Japan, 36 (1984), 295–313.

[13] Stepanov S. E., Shandra I. G., Geometry of infinitesimal harmonic transformations, Ann. Glob. Anal. Geom., 24 (2003), 291-299.

[14] Bourguignon J.-P., Ebin D. G., Marsden J. E., Sur le noyau des operaterus pseudo-diffirentiels à symbole surjectif et non injectif, C. R. Acad. Sc. Paris, Série A, 281 (1976), 867-870.

[15] Gicquaud R., Ngo Q. A., A new point of view on the solutions to the Einstein constraint equations with arbitrary mean curvature and small $TT$-tensor, Class. Quant. Grav. 31(19) (2014), 195014.

[16] Beig, R., Chruściel, P. T., All transverse and $TT$-tensors in flat spaces of any dimension, General Relativity and Gravitation, 50(5), (2018), 52.

[17] Rovenski V., Stepanov S., Tsyganok I., The Sampson Laplacian on negatively pinched Riemannian manifolds, International Electronic Journal of Geometry, 14(1), (2021), 91–99.

[18] Cao, X., Gursky, M. J., Tran, H., Curvature of the second kind and a conjecture of Nishikawa. Inventiones mathematicae, 233(3), (2023), 1347-1376.